\documentclass[runningheads]{llncs}
\usepackage[utf8]{inputenc}
\usepackage{amsmath}
\usepackage{graphicx}
\usepackage{algorithm}
\graphicspath{ {.} }
\usepackage{subfig}
\usepackage{svg}
\usepackage{amsfonts}
\usepackage{lineno}

\begin{document}
\title{On the Simple Quasi Crossing Number of $K_{11}$}

\author{Arjun Pitchanathan\inst{1} \and
Saswata Shannigrahi\inst{2}}

\authorrunning{A. Pitchanathan and S. Shannigrahi}

\institute{International Institute of Information Technology, Hyderabad 500032, India. \email{arjun.p@research.iiit.ac.in} \and
Indian Institute of Technology, Ropar 140001, India.
\email{saswata.shannigrahi@gmail.com}}

\maketitle

\begin{abstract}
\centering
We show that the simple quasi crossing number of $K_{11}$ is 4.
\end{abstract}

A {\it quasi-planar} graph \cite{quasiplanar-definition} is one that can be drawn in the plane without any triples of pairwise crossing edges. A drawing of a graph in the plane is {\it simple} if every pair of edges meets at most once, either at an intersection point or at a common endpoint. Accordingly, a graph is {\it simple quasi-planar} if it has a simple drawing in the plane without any triples of pairwise crossing edges. We define the {\it simple quasi crossing number} of a graph $G$, denoted by $cr_3(G)$, to be the minimum number of such triples in a simple drawing of $G$ in the plane.

It has been shown \cite{AT07} that a simple quasi-planar graph with $n \ge 4$ vertices has at most $6.5n - 20$ edges and this bound is tight up to a constant. We use this bound and follow the proof of the crossing number inequality \cite{cross_1,cross_2} to obtain a lower bound on the simple quasi crossing number of a graph as follows. 

Let $G$ be a graph with $n \ge 4$ vertices and $e$ edges. Consider a simple drawing of $G$ with $cr_3(G)$ triples of pairwise crossing edges. We can remove each such triple by removing an edge. In this way, we can obtain a simple quasi-planar graph with at least $e - cr_3(G)$ edges and $n$ vertices. By the above-mentioned bound, we have
\begin{equation}
    cr_3(G) \ge e - 6.5n + 20.
\end{equation}

We improve this bound by the probabilistic method, as in the proof of the crossing number inequality. Let $p$ be a parameter between 0 and 1, to be chosen later. Consider a random subgraph of $H$ obtained by including each vertex of $G$ independently with a probability $p$, and including the edges for which both vertices are included. Let $n_H$, $e_H$ and $cr_3(H)$ be the random variables denoting the number of vertices, number of edges and the simple quasi crossing number of $H$, respectively. By applying the inequality (1) and taking expectations, we obtain $\mathbb{E}[cr_3(H)] \ge \mathbb{E}[e_H] - 6.5 \  \mathbb{E}[n_H] + 20$. By the independence of the choices, we have $E[e_H] = p^2 e$ and $\mathbb{E}[n_H] = pn$. In any triple of pairwise crossing edges, there are exactly six distinct vertices involved. Therefore, we have $\mathbb{E}[cr_3(H)] \le p^6 cr_3(G)$. Setting $p = \alpha n / e$ and simplifying, we obtain
\begin{equation}
    cr_3(G) \ge \bigg( \frac{\alpha - 6.5}{\alpha^5} \bigg) \frac{e^5}{n^4} + \frac{20 e^6}{\alpha^6 n^6}
\end{equation}
for graphs satisfying $e \ge \alpha n$ since the probability $p$ must be at most 1. The value of $\alpha$ maximizing $(\alpha - 6.5)/\alpha^5$ is 8.125, which implies that $cr_3(G) \ge (1.625/8.125^5) e^5/n^4 + (20/8.125^6) e^6 / n^6$ for graphs satisfying $e \ge 8.125n$.

In this paper, we are particularly interested in the values of $cr_3(K_n)$ where $K_n$ denotes the complete graph on $n$ vertices. For $n \le 10$, $K_n$ is known to be simple quasi-planar \cite{quasiplanar,k10} and so we have $cr_3(K_n) = 0$. For $n = 11$, we obtain $cr_3(K_{11}) \ge 3.5$ from (1). Therefore, there must be at least four triples of pairwise crossing edges in any simple drawing of $K_{11}$. In the following, we present a drawing (Figure d) that shows that $cr_3(K_{11}) = 4$. In each of the Figures a-d, the triples of pairwise crossing edges introduced in the figure are marked with red circles.

\begin{figure}[h]
    \centering
    
    \subfloat[With the initial edges]{
        \def\svgwidth{0.41\textwidth}
    	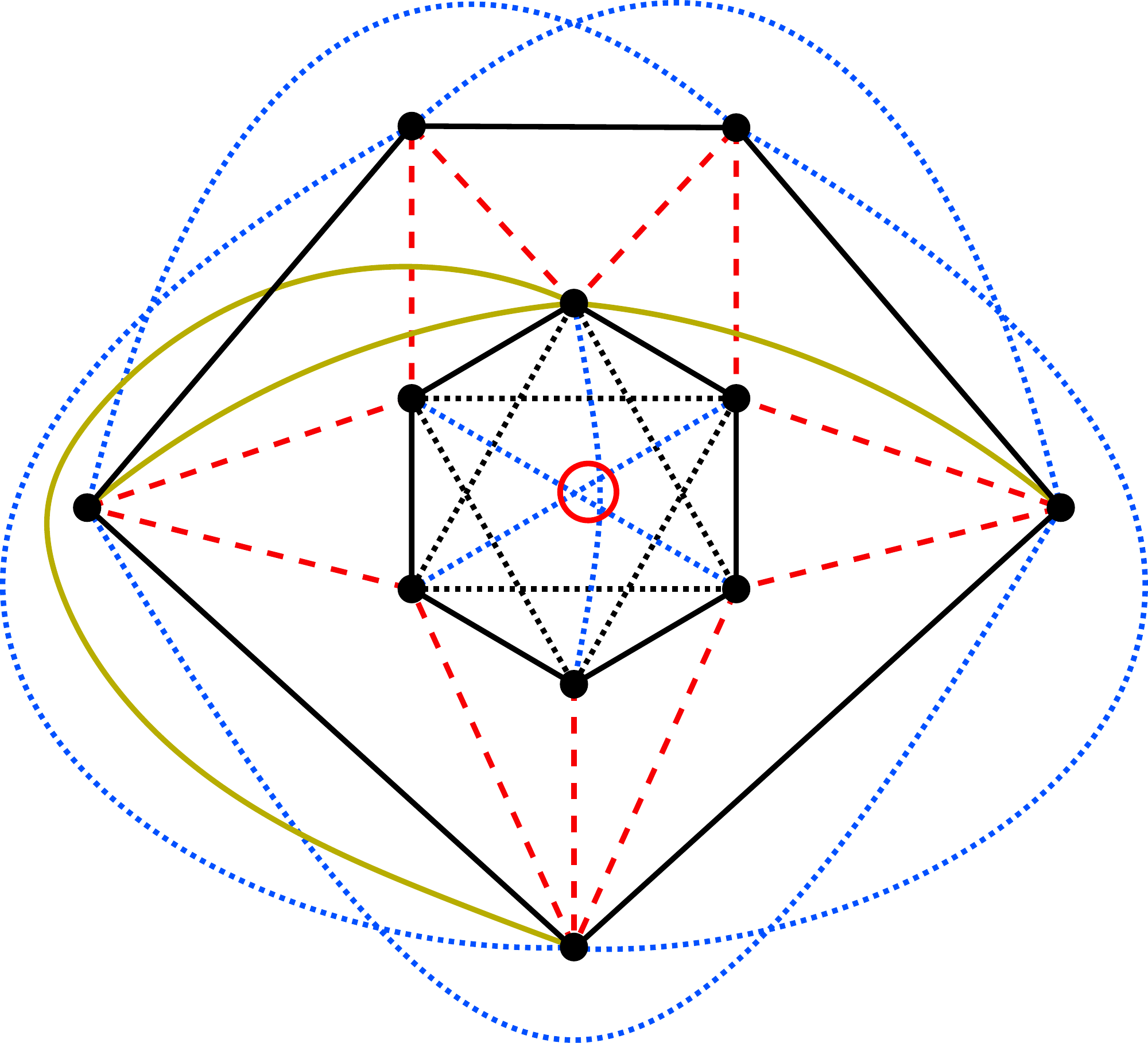
    }
    \qquad
    \subfloat[After adding pink edges]{
        \def\svgwidth{0.41\textwidth}
    	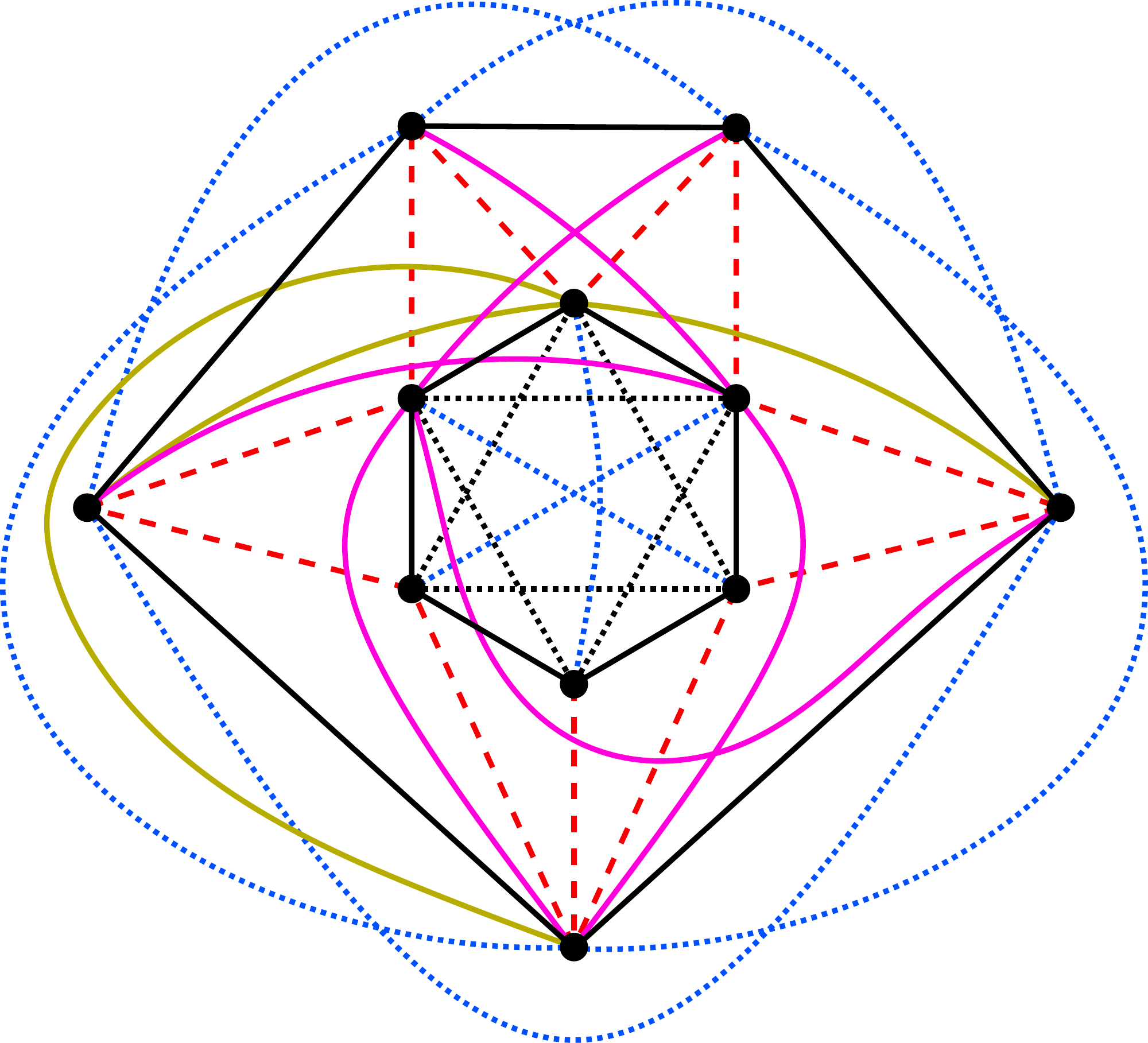
    }
    
    \subfloat[After adding dark blue edges]{
        \def\svgwidth{0.41\textwidth}
    	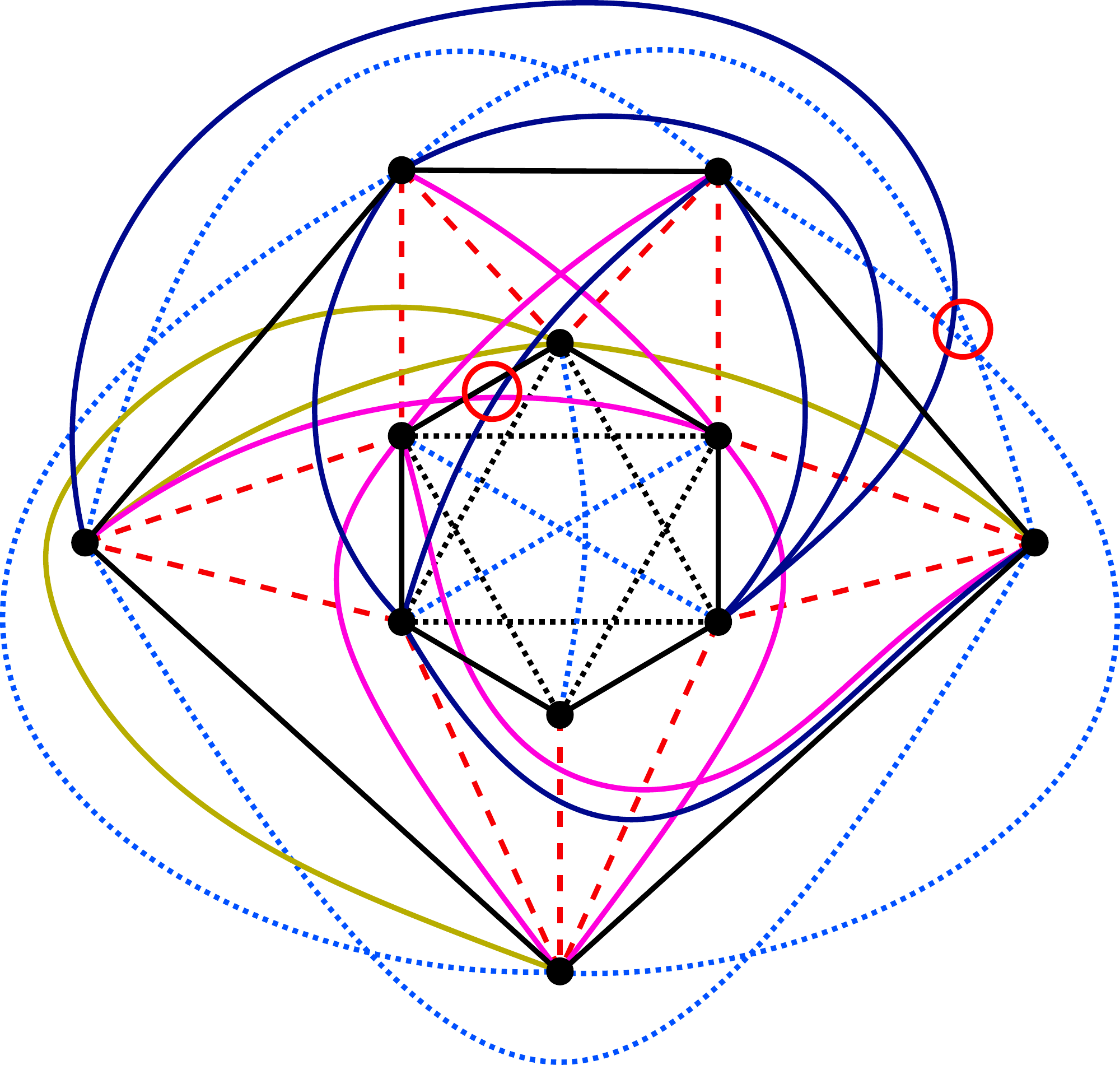
    }
    \qquad
    \subfloat[The complete drawing of $K_{11}$]{
        \def\svgwidth{0.41\textwidth}
    	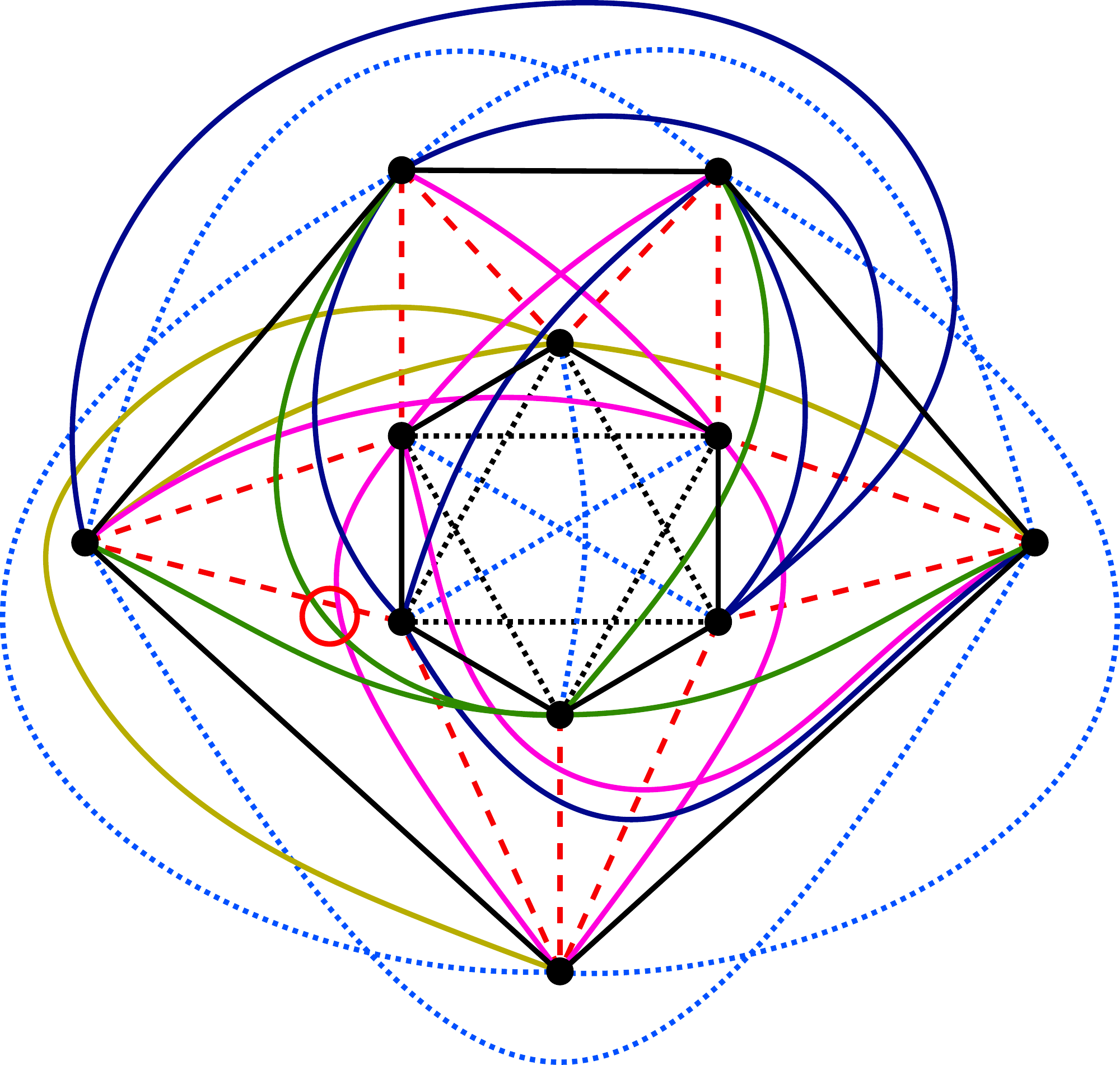
    }
\end{figure}

For $n \ge 12$, we are not aware of the exact values of $cr_3(K_n)$. In this case, the best known lower bounds can be obtained from the inequalities (1) and (2). It is an open problem to obtain non-trivial upper bounds on $cr_3(K_n)$ for $n \ge 12$. Another open problem is to find a general drawing that provides a good upper bound on $cr_3(K_n)$ for large $n$.



\end{document}